\documentclass[12pt,a4paper]{article}

\usepackage{graphicx}

\usepackage{amsmath}

\usepackage{amssymb}
\usepackage{hyperref}

\usepackage{amsthm}
\usepackage{geometry}
\usepackage{verbatim}
\usepackage{float}

\usepackage{setspace}
\spacing{1.25}
\usepackage{fancyhdr}

\usepackage{color} 

\graphicspath{{pics/}}

\newcommand{\al}{\alpha}

\newcommand{\vphi}{\varphi}

\newcommand{\be}{\beta}

\newcommand{\ga}{\gamma}

\newcommand{\de}{\delta}

\newcommand{\om}{\omega}

\newcommand{\na}{\nabla}

\newcommand{\NA}{\nabla}

\newcommand{\bs}{\boldsymbol}

\newcommand{\ra}{\rightarrow}

\newcommand{\lra}{\longrightarrow}

\newcommand{\Ra}{\Rightarrow}

\newcommand{\xra}{\xrightarrow}

\newcommand{\xlra}{\xlongrightarrow}

\newcommand{\rgl}{\rangle}

\newcommand{\lgl}{\langle}

\newcommand{\dash}{\textrm{-}}

\newcommand{\ot}{\otimes}

\newcommand{\bpf}{\begin{proof}}

\newcommand{\epf}{\end{proof}}

\newcommand{\bthm}{\begin{thm}}

\newcommand{\ethm}{\end{thm}}

\newcommand{\bprop}{\begin{prop}}

\newcommand{\eprop}{\end{prop}}

\newcommand{\bcor}{\begin{cor}}

\newcommand{\ecor}{\end{cor}}

\newcommand{\blem}{\begin{lem}}

\newcommand{\elem}{\end{lem}}

\newcommand{\bdefn}{\begin{defn}}

\newcommand{\edefn}{\end{defn}}

\newcommand{\bexmp}{\begin{exmp}}

\newcommand{\eexmp}{\end{exmp}}

\newcommand{\brem}{\begin{rem}}

\newcommand{\erem}{\end{rem}}

\newcommand{\bdia}{\begin{displaymath}\xymatrix}

\newcommand{\edia}{\end{displaymath}}

\newcommand{\beq}{\begin{equation*}\begin{aligned}}

\newcommand{\eeq}{\end{aligned}\end{equation*}}

\newcommand{\bref}{\textbf{Ref}}

\newcommand{\intg}{\mathbb{Z}}

\newcommand{\real}{\mathbb{R}}

\newcommand{\comp}{\mathbb{C}}

\newcommand{\quot}{\mathbb{H}}

\newcommand{\afv}{\mathbb{A}}

\newcommand{\prv}{\mathbb{P}}

\newcommand{\mco}{\mathcal{O}}

\newcommand{\mcc}{\mathcal{C}}

\newcommand{\mcf}{\mathcal{F}}

\newcommand{\mcg}{\mathcal{G}}

\newcommand{\mcs}{\mathcal{S}}

\newcommand{\cp}{\mathbb{CP}}

\newcommand{\mfo}{\mathfrak{O}}

\newcommand{\mfg}{\mathfrak{g}}

\newcommand{\msa}{\mathscr{A}}

\newcommand{\msr}{\mathscr{R}}

\newcommand{\msg}{\mathscr{G}}

\newcommand{\msd}{\mathscr{D}}

\newcommand{\itbf}{\item\textbf}

\newcommand{\seqa}{a_1,...,a_}

\newcommand{\seqx}{x_1,...,x_}

\newcommand{\seqy}{y_1,...,y_}

\newcommand{\seqf}{f_1,...,f_}

\newcommand{\cred}{\textcolor{red}}

\newcommand{\cblue}{\textcolor{blue}}

\newcommand{\mfa}{\mathfrak{a}}

\newcommand{\mfb}{\mathfrak{b}}

\newcommand{\mfm}{\mathfrak{m}}

\newcommand{\mfn}{\mathfrak{n}}

\newcommand{\mfp}{\mathfrak{p}}

\newcommand{\Af}{A_{(f)}}

\newtheorem{thm}{\textbf {Theorem}}[section]

\newtheorem{cor}[thm]{\textbf{Corollary}}

\newtheorem{prop}[thm]{\textbf{Proposition}}

\newtheorem{lem}[thm]{\textbf{Lemma}}

\theoremstyle{definition}

\newtheorem{defn}[thm]{\textbf{Definition}}

\newtheorem{exmp}[thm]{Example}

\newtheorem{notn}[thm]{Notation}

\newtheorem{conv}[thm]{Convention}

\theoremstyle{remark}

\newtheorem{rem}[thm]{Remark}

\def\cok{\operatorname{Coker}}

\newcommand{\txi}{\tilde{\xi}}

\newcommand{\bxi}{\bar{\xi}}

\newcommand{\bz}{\bar{z}}

\DeclareMathOperator{\tr}{trunk}

\title{Trunk of Satellite and Companion Knots}

\author{Nithin Kavi, Wendy Wu, and Zhenkun Li}

\date{}

\def\allfiles{}
\begin{document}

\bibliographystyle{plain}
\maketitle

\begin{abstract}
We study the knot invariant called trunk, as defined by Ozawa \cite{ozawa2010waist}, and the relation of the trunk of a satellite knot with the trunk of its companion knot. Our first result is $\tr(K) \geq n \cdot \tr(J)$ where $\tr(\cdot)$ denotes trunk of a knot, $K$ is a satellite knot with companion $J$, and $n$ is the winding number of $K$. To upgrade winding number to wrapping number, denoted by $m$, an extra factor of $\frac{1}{2}$ is necessary in our second result $\tr(K) > \frac{1}{2} m\cdot \tr(J)$ as $m \geq n$. We also discuss generalizations of the second result.

\end{abstract}



\ifx\allfiles\undefined

\documentclass[12pt,a4paper]{article}


\usepackage{graphicx}

\usepackage{amsmath}

\usepackage{amssymb}

\usepackage{amsthm}

\usepackage{geometry}

\DeclareMathOperator{\tr}{trunk}

\usepackage{fancyhdr}

\usepackage{color} 









\newcommand{\al}{\alpha}

\newcommand{\vphi}{\varphi}

\newcommand{\be}{\beta}

\newcommand{\ga}{\gamma}

\newcommand{\de}{\delta}

\newcommand{\om}{\omega}

\newcommand{\na}{\nabla}

\newcommand{\NA}{\nabla}

\newcommand{\bs}{\boldsymbol}

\newcommand{\ra}{\rightarrow}

\newcommand{\lra}{\longrightarrow}

\newcommand{\Ra}{\Rightarrow}

\newcommand{\xra}{\xrightarrow}

\newcommand{\xlra}{\xlongrightarrow}

\newcommand{\rgl}{\rangle}

\newcommand{\lgl}{\langle}

\newcommand{\dash}{\textrm{-}}

\newcommand{\ot}{\otimes}

\newcommand{\bpf}{\begin{proof}}

\newcommand{\epf}{\end{proof}}

\newcommand{\bthm}{\begin{thm}}

\newcommand{\ethm}{\end{thm}}

\newcommand{\bprop}{\begin{prop}}

\newcommand{\eprop}{\end{prop}}

\newcommand{\bcor}{\begin{cor}}

\newcommand{\ecor}{\end{cor}}

\newcommand{\blem}{\begin{lem}}

\newcommand{\elem}{\end{lem}}

\newcommand{\bdefn}{\begin{defn}}

\newcommand{\edefn}{\end{defn}}

\newcommand{\bexmp}{\begin{exmp}}

\newcommand{\eexmp}{\end{exmp}}

\newcommand{\brem}{\begin{rem}}

\newcommand{\erem}{\end{rem}}

\newcommand{\bdia}{\begin{displaymath}\xymatrix}

\newcommand{\edia}{\end{displaymath}}

\newcommand{\beq}{\begin{equation*}\begin{aligned}}

\newcommand{\eeq}{\end{aligned}\end{equation*}}

\newcommand{\bref}{\textbf{Ref}}

\newcommand{\intg}{\mathbb{Z}}

\newcommand{\real}{\mathbb{R}}

\newcommand{\comp}{\mathbb{C}}

\newcommand{\quot}{\mathbb{H}}

\newcommand{\afv}{\mathbb{A}}

\newcommand{\prv}{\mathbb{P}}

\newcommand{\mco}{\mathcal{O}}

\newcommand{\mcc}{\mathcal{C}}

\newcommand{\mcf}{\mathcal{F}}

\newcommand{\mcg}{\mathcal{G}}

\newcommand{\mcs}{\mathcal{S}}

\newcommand{\cp}{\mathbb{CP}}

\newcommand{\mfo}{\mathfrak{O}}

\newcommand{\mfg}{\mathfrak{g}}

\newcommand{\msa}{\mathscr{A}}

\newcommand{\msr}{\mathscr{R}}

\newcommand{\msg}{\mathscr{G}}

\newcommand{\msd}{\mathscr{D}}

\newcommand{\itbf}{\item\textbf}

\newcommand{\seqa}{a_1,...,a_}

\newcommand{\seqx}{x_1,...,x_}

\newcommand{\seqy}{y_1,...,y_}

\newcommand{\seqf}{f_1,...,f_}

\newcommand{\cred}{\textcolor{red}}

\newcommand{\cblue}{\textcolor{blue}}

\newcommand{\mfa}{\mathfrak{a}}

\newcommand{\mfb}{\mathfrak{b}}

\newcommand{\mfm}{\mathfrak{m}}

\newcommand{\mfn}{\mathfrak{n}}

\newcommand{\mfp}{\mathfrak{p}}

\newcommand{\Af}{A_{(f)}}


\newtheorem{thm}{\textbf {Theorem}}[section]

\newtheorem{cor}[thm]{\textbf{Corollary}}

\newtheorem{prop}[thm]{\textbf{Proposition}}

\newtheorem{lem}[thm]{\textbf{Lemma}}

\newtheorem{conj}[thm]{Conjecture}

\newtheorem{prob}[thm]{Problem}

\newtheorem{exer}[thm]{Exercise}

\newtheorem{quest}[thm]{Question}

\theoremstyle{definition}

\newtheorem{defn}[thm]{\textbf{Definition}}

\newtheorem{defns}[thm]{Definitions}

\newtheorem{exmp}[thm]{Example}

\newtheorem{exmps}[thm]{Examples}

\newtheorem{var}[thm]{Variant}

\newtheorem{vars}[thm]{Variants}

\newtheorem{con}[thm]{Construction}

\newtheorem{notn}[thm]{Notation}

\newtheorem{notns}[thm]{Notations}

\newtheorem{conv}[thm]{Convention}

\theoremstyle{remark}

\newtheorem{rem}[thm]{Remark}

\newtheorem{rems}[thm]{Remarks}

\newtheorem{warn}[thm]{Warning}

\newtheorem{sch}[thm]{Scholium}

\newtheorem{expl}[thm]{Explanations}

\newtheorem*{theorem}{\textbf{Theorem}}

\newtheorem*{corollary}{\textbf{Corollary}}

\newtheorem*{proposition}{\textbf{Proposition}}

\newtheorem*{lemma}{\textbf{Lemma}}

\newtheorem*{example}{\textbf{Example}}

\def\cok{\operatorname{Coker}}

\newcommand{\txi}{\tilde{\xi}}

\newcommand{\bxi}{\bar{\xi}}

\newcommand{\bz}{\bar{z}}



\begin{document}

\bibliographystyle{plain}

\maketitle

\else

\fi




\section{Introduction}
Knots and links are core objects in the study of 3-manifolds. The most important tools to study them are their numerical and homological invariants. There is an important family of invariants for knots (and links) which come from Morse theory. Among them are bridge number, width and trunk. For these Morse-type invariants, an important question is how they behave under the operations of connected sum and taking satellites, which are ways to construct complicated knots out of simple ones. Understanding the behaviors of the invariants under those operations would then contribute to the study of the properties of knots and links.

In this paper, we use $k$ to denote a specific knot, $\hat{V}$ to mean an unknotted solid torus, $\hat{j}$ to mean $S^1 \times 0$, and $K$ to refer to a knot class. 

With this, we can define the satellite knot:

{\bf Definition 1.}
{\it Let $\hat{k}$ be a knot inside $\hat{V}$ such that $\hat{k}$ intersects any meridian disk in $\hat{V}$. Let $f$ be a smooth embedding from $\hat{V}$ to $S^3$ and let $f(\hat{j}) = j$ and $f(\hat{k}) = k$. The knot $k$ is a {\it satellite knot} with {\it companion} $j$.}

Bridge number was first introduced by Schubert \cite{schubert1954uber} in the 1950s, and it has broad connections and applications in many aspects of knot theory. Its behavior under satellite operations has been studied by Schubert \cite{schubert1954uber} and Schultens \cite{schultens2003additivity}:
$$b(K_1\sharp K_2)=b(K_1)+b(K_2)-1,$$
$$b(K)\geq m\cdot b(J).$$
Here $b(\cdot)$ is the bridge number of a knot and $\sharp$ means the connected sum, defined in \cite{rolfsen1976knots}. $K_1$ and $K_2$ are knot classes. In the second inequality, $K$ is a satellite knot with companion $J$ and wrapping number $m$.

Width was first defined by Gabai \cite{gabai1987foliations} in his proof of the Property R conjecture. It is closely related to the study of meridional surfaces in the knot complements by Thompson \cite{thompson1994position} and Wu \cite{wu2006position}. It was also an essential ingredient of the proof of the knot complement conjecture by Gordon and Luecke \cite{gordon1989knots}. Its behavior under the connected sum was studied by Blair and Tomova \cite{blair2013width}, Rieck and Sedgwick \cite{rieck2002thin} and Scharlemann and Schultens \cite{scharlemann20063}, among others:
$$\max\{\omega(K_1),\omega(K_2)\}\leq \omega(K_1\sharp K_2)\leq \omega(K_1)+\omega(K_2)-2.$$

However, the behavior of width under taking satellites still remains a mystery. In \cite{guo2018width}, Guo and Li proved one of Zupan's conjectures from \cite{alex2010lower}:
\begin{equation}\label{eq_constant_one}
    \omega(K)\geq n^2\cdot \omega(J),
\end{equation}
where $K$ is a satellite knot with companion $J$ and $n$ is the winding number of $K$. This is not fully satisfactory as there are many important examples including Whitehead doubles which all have winding number zero, so inequality (\ref{eq_constant_one}) has only trivial results. As the wrapping number is always non-zero, we expect to replace the winding number $n$ by the wrapping number $m$ in the inequality (\ref{eq_constant_one}), leading to the conjecture:

{\bf Conjecture 1.} (Zupan) {\it Suppose $K$ is a satellite knot with companion $J$ and wrapping number $m$. Then we have}
$$\omega(K)\geq m^2\cdot \omega(J).$$

The special case where $K$ is the Whitehead double was proved by Guo and Li \cite{li2018width} but the general case is still open.

In this paper we present our results on trunk, which can be regarded as a simplified version of width, and could help with studying width. We ﬁrst adapt the methods used to prove the main result in \cite{guo2018width} and apply them to the trunk of a knot to prove the following theorem:

{\bf Theorem 1.} {\it Suppose $K$ is a satellite knot with companion $J$ and winding number $n$. Then we have}
\begin{equation}\label{eq_trunk_constant_one}
    {\rm trunk}(K)\geq n\cdot{\rm trunk}(J).
\end{equation}

We also study the case for wrapping number and obtain a lower bound of $\tr(K)$ in terms of $\tr(J)$ and the wrapping number $m$. By definition, the wrapping number is the least geometric intersection number of $K$ with any meridian disk of the tubular neighborhood of $J$ and is always non-zero (by the definition of taking satellites). However, we do not get a result as strong as inequality (\ref{eq_trunk_constant_one}) when we use the wrapping number as we have a factor of a half in our bound:

{\bf Theorem 2.} {\it Suppose $K$ is a satellite knot with companion $J$ and wrapping number $m$. Then we have}
$${\rm trunk}(K) > \frac{1}{2}\cdot m\cdot{\rm trunk}(J).$$

We still make the following conjecture:

{\bf Conjecture 2.}
{\it Suppose $K$ is a satellite knot with companion $J$ and wrapping number $m$. Then we have}
$${\rm trunk}(K)\geq m\cdot{\rm trunk}(J).$$

To bound the trunk of $K$, we need to study the intersection of a particular knot $k$ with the regular level $h^{-1}(r)$ of the standard Morse function $h$ on $S^3$ and a regular value $r\in \real$. Since our knot $k$ is contained in a tubular neighborhood $V$ of the companion knot, we can first study the intersection $V\cap h^{-1}(r)$. By definition, if a connected component $P\subset V\cap h^{-1}(r)$ is a meridian disk, then it intersects $K$ at least $m$ times. Hence the proof of Theorem 2 relies on the following key lemma:

{\bf Key lemma.} {\it Among all the relevant components (defined more precisely in Section 4) of $V\cap h^{-1}(r)$, at least half of them are meridian disks.}

There are some topological requirements for the intersection $V\cap h^{-1}(r)$. These requirements tell us how the components of $V\cap h^{-1}(r)$ are arranged on the regular level $h^{-1}(r)$ which is a $2$-sphere. Then we translate this problem into a purely combinatorial one about arranging pieces on a $2$-sphere and prove the key lemma in that setting.

The paper is organized as follows: In Section 2 we introduce some basic definitions about knot invariants and satellite knots. In Section 3 we summarize the result in \cite{guo2018width} and prove $\tr(K) \geq n \cdot \tr(J)$. In Section 4 we explain how to translate the problem into combinatorics and prove the Key lemma. We also discuss the wrapping number further and make some slight generalizations of Theorem 2.

\section*{Acknowledgments}

The first two authors thank their mentor Zhenkun Li for offering guidance throughout this project. Additionally, they thank the MIT math department and the MIT PRIMES program for giving them the opportunity to conduct this research. The authors also thank the anonymous referee for providing feedback on this paper.



\ifx\allfiles\undefined

\bibliography{Index}

\end{document}

\fi


\ifx\allfiles\undefined

\documentclass[12pt,a4paper]{article}


\usepackage{graphicx}

\usepackage{amsmath}

\usepackage{amssymb}
\DeclareMathOperator{\tr}{trunk}
\usepackage{amsthm}

\usepackage{geometry}

\usepackage{fancyhdr}

\usepackage{color} 









\newcommand{\al}{\alpha}

\newcommand{\vphi}{\varphi}

\newcommand{\be}{\beta}

\newcommand{\ga}{\gamma}

\newcommand{\de}{\delta}

\newcommand{\om}{\omega}

\newcommand{\na}{\nabla}

\newcommand{\NA}{\nabla}

\newcommand{\bs}{\boldsymbol}

\newcommand{\ra}{\rightarrow}

\newcommand{\lra}{\longrightarrow}

\newcommand{\Ra}{\Rightarrow}

\newcommand{\xra}{\xrightarrow}

\newcommand{\xlra}{\xlongrightarrow}

\newcommand{\rgl}{\rangle}

\newcommand{\lgl}{\langle}

\newcommand{\dash}{\textrm{-}}

\newcommand{\ot}{\otimes}

\newcommand{\bpf}{\begin{proof}}

\newcommand{\epf}{\end{proof}}

\newcommand{\bthm}{\begin{thm}}

\newcommand{\ethm}{\end{thm}}

\newcommand{\bprop}{\begin{prop}}

\newcommand{\eprop}{\end{prop}}

\newcommand{\bcor}{\begin{cor}}

\newcommand{\ecor}{\end{cor}}

\newcommand{\blem}{\begin{lem}}

\newcommand{\elem}{\end{lem}}

\newcommand{\bdefn}{\begin{defn}}

\newcommand{\edefn}{\end{defn}}

\newcommand{\bexmp}{\begin{exmp}}

\newcommand{\eexmp}{\end{exmp}}

\newcommand{\brem}{\begin{rem}}

\newcommand{\erem}{\end{rem}}

\newcommand{\bdia}{\begin{displaymath}\xymatrix}

\newcommand{\edia}{\end{displaymath}}

\newcommand{\beq}{\begin{equation*}\begin{aligned}}

\newcommand{\eeq}{\end{aligned}\end{equation*}}

\newcommand{\bref}{\textbf{Ref}}

\newcommand{\intg}{\mathbb{Z}}

\newcommand{\real}{\mathbb{R}}

\newcommand{\comp}{\mathbb{C}}

\newcommand{\quot}{\mathbb{H}}

\newcommand{\afv}{\mathbb{A}}

\newcommand{\prv}{\mathbb{P}}

\newcommand{\mco}{\mathcal{O}}

\newcommand{\mcc}{\mathcal{C}}

\newcommand{\mcf}{\mathcal{F}}

\newcommand{\mcg}{\mathcal{G}}

\newcommand{\mcs}{\mathcal{S}}

\newcommand{\cp}{\mathbb{CP}}

\newcommand{\mfo}{\mathfrak{O}}

\newcommand{\mfg}{\mathfrak{g}}

\newcommand{\msa}{\mathscr{A}}

\newcommand{\msr}{\mathscr{R}}

\newcommand{\msg}{\mathscr{G}}

\newcommand{\msd}{\mathscr{D}}

\newcommand{\itbf}{\item\textbf}

\newcommand{\seqa}{a_1,...,a_}

\newcommand{\seqx}{x_1,...,x_}

\newcommand{\seqy}{y_1,...,y_}

\newcommand{\seqf}{f_1,...,f_}

\newcommand{\cred}{\textcolor{red}}

\newcommand{\cblue}{\textcolor{blue}}

\newcommand{\mfa}{\mathfrak{a}}

\newcommand{\mfb}{\mathfrak{b}}

\newcommand{\mfm}{\mathfrak{m}}

\newcommand{\mfn}{\mathfrak{n}}

\newcommand{\mfp}{\mathfrak{p}}

\newcommand{\Af}{A_{(f)}}


\newtheorem{thm}{\textbf {Theorem}}[section]

\newtheorem{cor}[thm]{\textbf{Corollary}}

\newtheorem{prop}[thm]{\textbf{Proposition}}

\newtheorem{lem}[thm]{\textbf{Lemma}}

\newtheorem{conj}[thm]{Conjecture}

\newtheorem{conv}[thm]{Convention}

\newtheorem{prob}[thm]{Problem}

\newtheorem{exer}[thm]{Exercise}

\newtheorem{quest}[thm]{Question}

\theoremstyle{definition}

\newtheorem{defn}[thm]{\textbf{Definition}}

\newtheorem{defns}[thm]{Definitions}

\newtheorem{exmp}[thm]{Example}

\newtheorem{exmps}[thm]{Examples}

\newtheorem{var}[thm]{Variant}

\newtheorem{vars}[thm]{Variants}

\newtheorem{con}[thm]{Construction}

\newtheorem{notn}[thm]{Notation}

\newtheorem{notns}[thm]{Notations}

\theoremstyle{remark}

\newtheorem{rem}[thm]{Remark}

\newtheorem{rems}[thm]{Remarks}

\newtheorem{warn}[thm]{Warning}

\newtheorem{sch}[thm]{Scholium}

\newtheorem{expl}[thm]{Explanations}

\newtheorem*{theorem}{\textbf{Theorem}}

\newtheorem*{corollary}{\textbf{Corollary}}

\newtheorem*{proposition}{\textbf{Proposition}}

\newtheorem*{lemma}{\textbf{Lemma}}

\newtheorem*{example}{\textbf{Example}}

\def\cok{\operatorname{Coker}}

\newcommand{\txi}{\tilde{\xi}}

\newcommand{\bxi}{\bar{\xi}}

\newcommand{\bz}{\bar{z}}



\begin{document}

\bibliographystyle{plain}

\else

\fi

\section{Preliminaries}
We will start with some necessary definitions.

\begin{defn}\label{defn_knot}
A {\it knot} is a smooth embedding
$$k:S^1\hookrightarrow S^3$$
where $S^1$ is the unit circle in $\mathbb{R}^2$ and $S^3$ is the unit sphere in $\mathbb{R}^4$.
\end{defn}

\begin{defn}\label{defn_knot_class}
A {\it knot class} is a set of knots that are all isotopic to each other. See \cite{rolfsen1976knots} for the definition of an isotopy.
\end{defn}


We shall fix a Morse function throughout the paper. We consider $h: S^3 \rightarrow \mathbb{R}$ to be the standard height function $h(x, y, z, w) = w$ restricted to the unit sphere $S^3 \subset \mathbb{R}^4$. The pre-images of $\pm 1$ are denoted by $\pm \infty$.

\begin{defn}\label{defn_Morse}
With the above notation, a knot $k$ is called {\it Morse} if the following composition is a Morse function:
$$h\circ k:S^1\ra \real.$$
\end{defn}

\begin{conv} The following conventions will be used throughout the paper:

(1) We will only consider knots that are Morse and whose critical points are all at different levels.

(2) Knots are denoted by a lowercase letter like $k$, while knot classes are generally denoted by a capital letter like $K$.

(3) By a knot we can either mean the embedding $S^1\ra S^3$ or the image of the embedding. We do not distinguish between them. 
\end{conv}

\begin{notn}\label{defn_1}

Let $k$ be a knot in $S^3$. Denote the critical levels of $k$ by $c_i$, and pick regular levels $r_i$ between two consecutive critical levels $c_i$ and $c_{i + 1}$, so that:

$$c_1 < r_1 < c_2 < r_2 < \ldots < c_{s - 1} < r_{s - 1} < c_s.$$

 For each regular level $r_i$, we define $w_i = |h^{-1}(r_i) \cap k|$, that is, $w_i$ is the number of intersections of this regular level with $k$.
\end{notn}

\subsection{Trunk and Width of Knots}

Now we will define two invariants, {\it trunk} and {\it width}, for knots and knot classes.

\begin{defn}\mbox{}
\begin{itemize}
    \item The {\it trunk number} of a knot $k$ is given by $\tr(k) = \displaystyle{\max_{1 \leq i \leq s - 1} w_i(k)}$.
    \item The {\it width number} of a knot $k$ is given by $$\omega(k) = \sum_{ i = 1}^{s - 1} w_i.$$
\end{itemize}
\end{defn}



\begin{defn}\label{defn_2} We can extend the definitions of the trunk and width numbers of a knot to also apply to knot classes:
\begin{itemize}
    
    \item The {\it trunk number of a knot class} $K$ is given by $\tr(K) = \displaystyle{\min_{k \in K} \tr(k)}$.
    \item The {\it width number of a knot class} $K$ is given by $\omega(K) = \displaystyle{\min_{k \in K} \omega(k)}$.
\end{itemize}


\end{defn}

\begin{figure}[H]
    \centering
    \includegraphics[scale = .4]{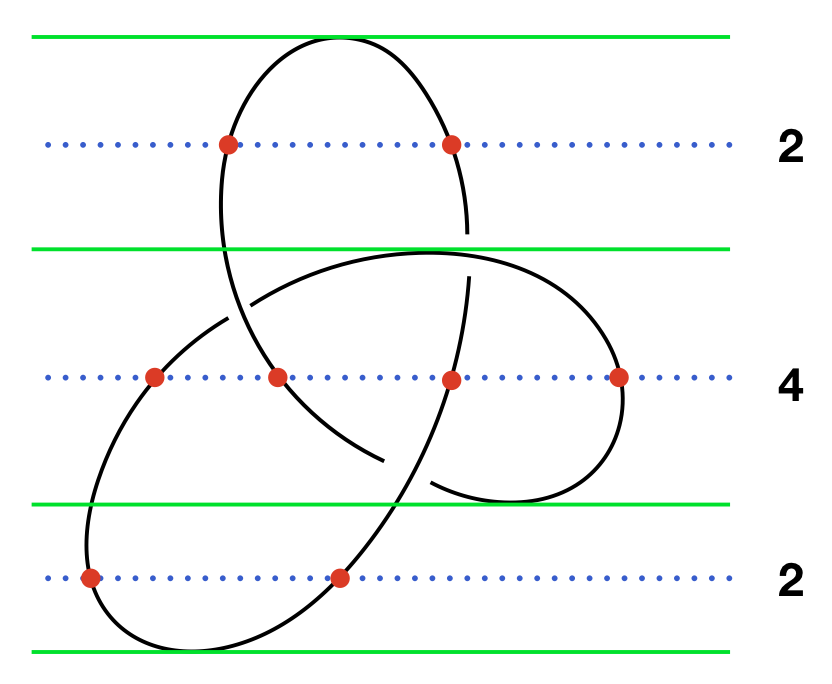}
    \caption{The width and trunk of a trefoil knot are 8 and 4, respectively.}
    \label{fig_trunk}
\end{figure}

\begin{exmp}\label{tref}
Suppose $K$ is the trefoil knot. Let $k$ be the particular embedding as in Figure \ref{fig_trunk}. 
Its width and trunk happen to be those of the trefoil knot class.
\end{exmp}

\begin{defn}\label{defn_essential} Let $T^2 = S^1 \times S^1$ be a two dimensional torus. A curve $\al \subseteq T^2$ is {\it inessential} if $\al = \partial D$ for some $D \subseteq T^2$ and is {\it essential} otherwise.

\end{defn}

\subsection{Satellite and Companion Knots}

Here, we will define the process of forming a {\it satellite knot}, which is one of the main ways to construct complicated knots from simple ones. We also define two invariants pertaining to satellite knots inside a solid torus.

\begin{defn}
Let $\hat{V} = S^1 \times D^2,$ and $\hat{j} = S^1 \times 0.$ A {\it meridian disk} of $\hat{V}$ is a properly embedded disk $D$ whose boundary $\partial{D}\subset\partial{\hat{V}}$ is essential on $\partial{\hat{V}}$.
\end{defn}

\begin{defn}\label{satellite}
Let $\hat{k}$ be a knot inside $\hat{V}$ such that $\hat{k}$ intersects any meridian disk in $\hat{V}$. Let $f$ be a smooth embedding from $\hat{V}$ to $S^3$ and let $f(\hat{j}) = j$ and $f(\hat{k}) = k$. The knot $k$ is a {\it satellite knot} with {\it companion} $j$.

\begin{figure}[H]
    \centering
    \includegraphics[width = 0.8\linewidth]{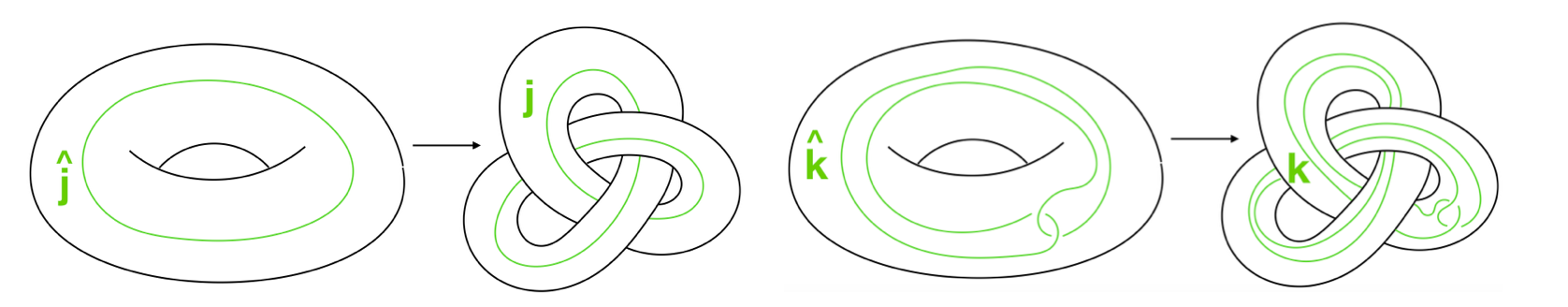}
    \caption{A satellite knot with companion being a trefoil.}
    \label{fig_satellite knot}
\end{figure}
\end{defn}

\begin{defn}\label{winding}

The {\it winding number} $n$ is the absolute value of the sum of all algebraic intersections of any fixed meridian disk with $\hat{k}$. The {\it wrapping number} $m$ is the minimal geometric intersection number of $k$ with any meridian disk. 

\end{defn}

\begin{figure}[h]
    \centering
    \includegraphics[scale = 0.3]{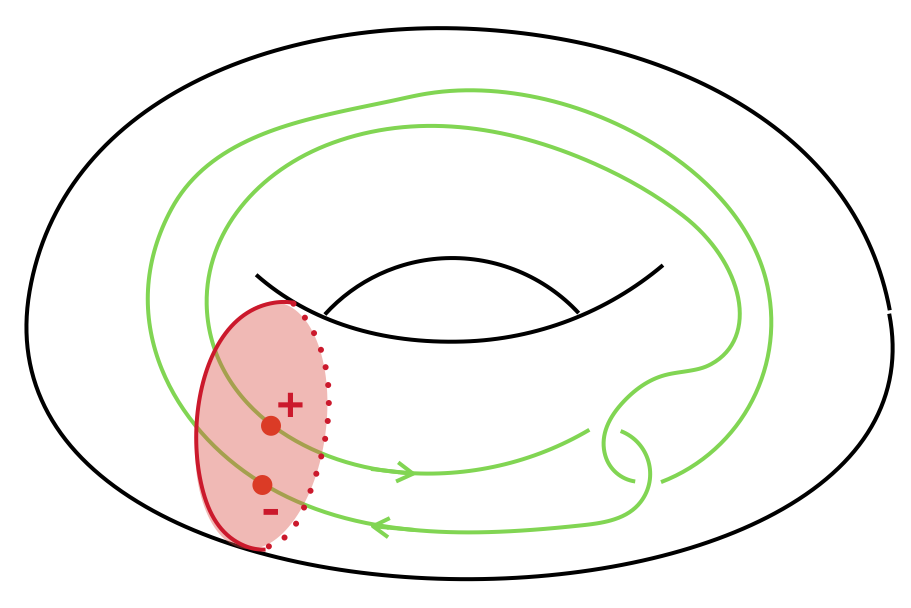}
    \caption{Satellite knot with winding number $0$ and wrapping number $2$.}
    \label{fig_winding_wrapping}
\end{figure}


\ifx\allfiles\undefined

\bibliography{Index}

\end{document}

\fi


\ifx\allfiles\undefined

\documentclass[12pt,a4paper]{article}


\usepackage{graphicx}

\usepackage{amsmath}

\usepackage{amssymb}

\usepackage{amsthm}

\usepackage{geometry}

\usepackage{fancyhdr}

\usepackage{color} 









\newcommand{\al}{\alpha}

\newcommand{\vphi}{\varphi}

\newcommand{\be}{\beta}

\newcommand{\ga}{\gamma}

\newcommand{\de}{\delta}

\newcommand{\om}{\omega}

\newcommand{\na}{\nabla}

\newcommand{\NA}{\nabla}

\newcommand{\bs}{\boldsymbol}

\newcommand{\ra}{\rightarrow}

\newcommand{\lra}{\longrightarrow}

\newcommand{\Ra}{\Rightarrow}

\newcommand{\xra}{\xrightarrow}

\newcommand{\xlra}{\xlongrightarrow}

\newcommand{\rgl}{\rangle}

\newcommand{\lgl}{\langle}

\newcommand{\dash}{\textrm{-}}

\newcommand{\ot}{\otimes}

\newcommand{\bpf}{\begin{proof}}

\newcommand{\epf}{\end{proof}}

\newcommand{\bthm}{\begin{thm}}

\newcommand{\ethm}{\end{thm}}

\newcommand{\bprop}{\begin{prop}}

\newcommand{\eprop}{\end{prop}}

\newcommand{\bcor}{\begin{cor}}

\newcommand{\ecor}{\end{cor}}

\newcommand{\blem}{\begin{lem}}

\newcommand{\elem}{\end{lem}}

\newcommand{\bdefn}{\begin{defn}}

\newcommand{\edefn}{\end{defn}}

\newcommand{\bexmp}{\begin{exmp}}

\newcommand{\eexmp}{\end{exmp}}

\newcommand{\brem}{\begin{rem}}

\newcommand{\erem}{\end{rem}}

\newcommand{\bdia}{\begin{displaymath}\xymatrix}

\newcommand{\edia}{\end{displaymath}}

\newcommand{\beq}{\begin{equation*}\begin{aligned}}

\newcommand{\eeq}{\end{aligned}\end{equation*}}

\newcommand{\bref}{\textbf{Ref}}

\newcommand{\intg}{\mathbb{Z}}

\newcommand{\real}{\mathbb{R}}

\newcommand{\comp}{\mathbb{C}}

\newcommand{\quot}{\mathbb{H}}

\DeclareMathOperator{\tr}{trunk}
\newcommand{\afv}{\mathbb{A}}

\newcommand{\prv}{\mathbb{P}}

\newcommand{\mco}{\mathcal{O}}

\newcommand{\mcc}{\mathcal{C}}

\newcommand{\mcf}{\mathcal{F}}

\newcommand{\mcg}{\mathcal{G}}

\newcommand{\mcs}{\mathcal{S}}

\newcommand{\cp}{\mathbb{CP}}

\newcommand{\mfo}{\mathfrak{O}}

\newcommand{\mfg}{\mathfrak{g}}

\newcommand{\msa}{\mathscr{A}}

\newcommand{\msr}{\mathscr{R}}

\newcommand{\msg}{\mathscr{G}}

\newcommand{\msd}{\mathscr{D}}

\newcommand{\itbf}{\item\textbf}

\newcommand{\seqa}{a_1,...,a_}

\newcommand{\seqx}{x_1,...,x_}

\newcommand{\seqy}{y_1,...,y_}

\newcommand{\seqf}{f_1,...,f_}

\newcommand{\cred}{\textcolor{red}}

\newcommand{\cblue}{\textcolor{blue}}

\newcommand{\mfa}{\mathfrak{a}}

\newcommand{\mfb}{\mathfrak{b}}

\newcommand{\mfm}{\mathfrak{m}}

\newcommand{\mfn}{\mathfrak{n}}

\newcommand{\mfp}{\mathfrak{p}}

\newcommand{\Af}{A_{(f)}}


\newtheorem{thm}{\textbf {Theorem}}[section]

\newtheorem{cor}[thm]{\textbf{Corollary}}

\newtheorem{prop}[thm]{\textbf{Proposition}}

\newtheorem{lem}[thm]{\textbf{Lemma}}

\newtheorem{conj}[thm]{Conjecture}

\newtheorem{prob}[thm]{Problem}

\newtheorem{exer}[thm]{Exercise}

\newtheorem{quest}[thm]{Question}

\theoremstyle{definition}

\newtheorem{defn}[thm]{\textbf{Definition}}

\newtheorem{defns}[thm]{Definitions}

\newtheorem{exmp}[thm]{Example}

\newtheorem{exmps}[thm]{Examples}

\newtheorem{var}[thm]{Variant}

\newtheorem{vars}[thm]{Variants}

\newtheorem{con}[thm]{Construction}

\newtheorem{notn}[thm]{Notation}

\newtheorem{notns}[thm]{Notations}

\newtheorem{conv}[thm]{Convention}

\theoremstyle{remark}

\newtheorem{rem}[thm]{Remark}

\newtheorem{rems}[thm]{Remarks}

\newtheorem{warn}[thm]{Warning}

\newtheorem{sch}[thm]{Scholium}

\newtheorem{expl}[thm]{Explanations}

\newtheorem*{theorem}{\textbf{Theorem}}

\newtheorem*{corollary}{\textbf{Corollary}}

\newtheorem*{proposition}{\textbf{Proposition}}

\newtheorem*{lemma}{\textbf{Lemma}}

\newtheorem*{example}{\textbf{Example}}

\def\cok{\operatorname{Coker}}

\newcommand{\txi}{\tilde{\xi}}

\newcommand{\bxi}{\bar{\xi}}

\newcommand{\bz}{\bar{z}}

\DeclareMathOperator{\tr}{tr}



\begin{document}

\bibliographystyle{plain}

\else

\fi

\section{Bounding Trunk with Winding Number}
In this section we will prove Theorem \ref{thm_bound_with_winding_number} as well as review the main ideas in \cite{guo2018width}.

\begin{thm}\label{thm_bound_with_winding_number}
Suppose $K$ is a satellite knot with a non-trivial companion $J$ and the winding number is $n$. Then we have: $\tr(K) \ge n \cdot \tr(J)$.
\end{thm}

\begin{proof}
We pick a knot $k \in K$ such that $\tr(k) = \tr(K)$. There will be a corresponding companion $j$ and a solid torus $V$ containing $j$ and $k$ as in Definition \ref{satellite}. As in \cite{guo2018width}, we can assume that $h|_{\partial{V}}$ is Morse and all critical points of $h|_{\partial{V}}$ are in distinct levels. From the Pop Over Lemma in \cite{schultens2003additivity}, we can also assume that $V$ does not contain the two critical points $\pm\infty$ of $S^3$. Let $c_1,...,c_s$ be all critical values of $h|_{\partial{V}}$. We define

$$M=V\backslash(\mathop{\bigcup}_{i=1}^s h^{-1}(c_i)).$$

We construct a graph out of the decomposition of $V$ induced by $M$ where vertices correspond to components of $M$ and two vertices are connected by an edge if the corresponding two components of $M$ are separated by a critical level. We call this graph $\Gamma_{R}(V)$. 

\brem
The graph of this type was first introduced in the paper \cite{scharlemann20063} by Scharlemann and Schultens. Later Guo and Li made a similar construction in Definition 2.4 of  \cite{guo2018width}. Here we use the same construction as in Guo and Li's paper, where this graph is called a Reeb graph.
\erem

\begin{prop}[Guo, Li \cite{guo2018width}]\label{prop_property_of_the_graph}
The graph $\Gamma_{R}(V)$ has the following properties:

(1) There is a unique loop $l\subset \Gamma_{R}(V)$. We can embed $l$ into $V$.

(2) The loop $l$ represents a generator in $H_1(V)\cong \intg$.

(3) The loop $l\subset V$ can be also considered as a knot $l\subset S^3$ and its knot class $L$ is a connected sum of the companion $J$ with another knot $J'$:
$$L=J\# J'.$$

These properties follow from Lemma 3.6, Proposition 3.7, and Lemma 4.3 in \cite{guo2018width}.

\end{prop}

\begin{thm}[Davies, Zupan \cite{davies2017natural}]\label{thm_zupan} For two knots $K_1, K_2$ we have $\tr(K_1 \# K_2) = \max \{\tr(K_1), \tr(K_2)\}$.

\end{thm}

From Proposition \ref{prop_property_of_the_graph} and Theorem \ref{thm_zupan}, we have $\tr(L) = \tr(J \# J') \geq \tr(J)$. Therefore, to prove Theorem \ref{thm_bound_with_winding_number}, we simply need to show that $\tr(K) \geq n \cdot \tr(L)$. We will need the following two lemmas, which follow from Lemma 4.4 in \cite{guo2018width}:

\begin{lem}[Guo, Li \cite{guo2018width}]\label{lem_li_2} We can isotope $l$ into such a position $l'$, also contained in $V$, so that for any regular value $r\in\real$, we have the following property: suppose all components of intersection $h^{-1}(r)\cap V$ are
$$h^{-1}(r)\cap V=P_1\cup P_2\cup ... P_{t}.$$
Then each component $P_i$ intersects $l$ at most once.

\end{lem}

\begin{lem}[Guo, Li \cite{guo2018width}]\label{lem_1i_1} Let $l'$ be given as in the above lemma. Given a planar surface $P$ where $|P \cap l'| = 1$, we have $|P \cap k| \geq n$.

\end{lem}

We can choose a regular level $r$ such that $|h^{-1}(r) \cap l| = \tr(l)$. The above two lemmas apply here to conclude:
$$|h^{-1}(r) \cap k| = n\cdot \tr(l).$$
\end{proof}

\ifx\allfiles\undefined

\newpage

\bibliography{Index}

\end{document}

\fi

\ifx\allfiles\undefined

\documentclass[12pt,a4paper]{article}


\usepackage{graphicx}

\usepackage{amsmath}

\usepackage{amssymb}

\usepackage{amsthm}

\usepackage{geometry}

\usepackage{fancyhdr}

\usepackage{color} 









\newcommand{\al}{\alpha}

\newcommand{\vphi}{\varphi}

\newcommand{\be}{\beta}

\newcommand{\ga}{\gamma}

\newcommand{\de}{\delta}

\newcommand{\om}{\omega}

\newcommand{\na}{\nabla}

\newcommand{\NA}{\nabla}

\newcommand{\bs}{\boldsymbol}

\newcommand{\ra}{\rightarrow}

\newcommand{\lra}{\longrightarrow}

\newcommand{\Ra}{\Rightarrow}

\newcommand{\xra}{\xrightarrow}

\newcommand{\xlra}{\xlongrightarrow}

\newcommand{\rgl}{\rangle}

\newcommand{\lgl}{\langle}

\newcommand{\dash}{\textrm{-}}

\newcommand{\ot}{\otimes}

\newcommand{\bpf}{\begin{proof}}

\newcommand{\epf}{\end{proof}}

\newcommand{\bthm}{\begin{thm}}

\newcommand{\ethm}{\end{thm}}

\newcommand{\bprop}{\begin{prop}}

\newcommand{\eprop}{\end{prop}}

\newcommand{\bcor}{\begin{cor}}

\newcommand{\ecor}{\end{cor}}

\newcommand{\blem}{\begin{lem}}

\newcommand{\elem}{\end{lem}}

\newcommand{\bdefn}{\begin{defn}}

\newcommand{\edefn}{\end{defn}}

\newcommand{\bexmp}{\begin{exmp}}

\newcommand{\eexmp}{\end{exmp}}

\newcommand{\brem}{\begin{rem}}

\newcommand{\erem}{\end{rem}}

\newcommand{\bdia}{\begin{displaymath}\xymatrix}

\newcommand{\edia}{\end{displaymath}}

\newcommand{\beq}{\begin{equation*}\begin{aligned}}

\newcommand{\eeq}{\end{aligned}\end{equation*}}

\newcommand{\bref}{\textbf{Ref}}

\newcommand{\intg}{\mathbb{Z}}

\newcommand{\real}{\mathbb{R}}

\newcommand{\comp}{\mathbb{C}}

\newcommand{\quot}{\mathbb{H}}

\newcommand{\afv}{\mathbb{A}}

\newcommand{\prv}{\mathbb{P}}

\newcommand{\mco}{\mathcal{O}}

\newcommand{\mcc}{\mathcal{C}}

\newcommand{\mcf}{\mathcal{F}}

\newcommand{\mcg}{\mathcal{G}}

\newcommand{\mcs}{\mathcal{S}}

\newcommand{\cp}{\mathbb{CP}}

\newcommand{\mfo}{\mathfrak{O}}

\newcommand{\mfg}{\mathfrak{g}}

\newcommand{\msa}{\mathscr{A}}

\newcommand{\msr}{\mathscr{R}}

\newcommand{\msg}{\mathscr{G}}

\newcommand{\msd}{\mathscr{D}}

\newcommand{\itbf}{\item\textbf}

\newcommand{\seqa}{a_1,...,a_}

\newcommand{\seqx}{x_1,...,x_}

\newcommand{\seqy}{y_1,...,y_}

\newcommand{\seqf}{f_1,...,f_}

\DeclareMathOperator{\tr}{trunk}
\newcommand{\cred}{\textcolor{red}}

\newcommand{\cblue}{\textcolor{blue}}

\newcommand{\mfa}{\mathfrak{a}}

\newcommand{\mfb}{\mathfrak{b}}

\newcommand{\mfm}{\mathfrak{m}}

\newcommand{\mfn}{\mathfrak{n}}

\newcommand{\mfp}{\mathfrak{p}}

\newcommand{\Af}{A_{(f)}}


\newtheorem{thm}{\textbf {Theorem}}[section]

\newtheorem{cor}[thm]{\textbf{Corollary}}

\newtheorem{prop}[thm]{\textbf{Proposition}}

\newtheorem{lem}[thm]{\textbf{Lemma}}

\newtheorem{conj}[thm]{Conjecture}

\newtheorem{prob}[thm]{Problem}

\newtheorem{exer}[thm]{Exercise}

\newtheorem{quest}[thm]{Question}

\theoremstyle{definition}

\newtheorem{defn}[thm]{\textbf{Definition}}

\newtheorem{defns}[thm]{Definitions}

\newtheorem{exmp}[thm]{Example}

\newtheorem{exmps}[thm]{Examples}

\newtheorem{var}[thm]{Variant}

\newtheorem{vars}[thm]{Variants}

\newtheorem{con}[thm]{Construction}

\newtheorem{notn}[thm]{Notation}

\newtheorem{notns}[thm]{Notations}

\newtheorem{conv}[thm]{Convention}

\theoremstyle{remark}

\newtheorem{rem}[thm]{Remark}

\newtheorem{rems}[thm]{Remarks}

\newtheorem{warn}[thm]{Warning}

\newtheorem{sch}[thm]{Scholium}

\newtheorem{expl}[thm]{Explanations}

\newtheorem*{theorem}{\textbf{Theorem}}

\newtheorem*{corollary}{\textbf{Corollary}}

\newtheorem*{proposition}{\textbf{Proposition}}

\newtheorem*{lemma}{\textbf{Lemma}}

\newtheorem*{example}{\textbf{Example}}

\def\cok{\operatorname{Coker}}

\newcommand{\txi}{\tilde{\xi}}

\newcommand{\bxi}{\bar{\xi}}

\newcommand{\bz}{\bar{z}}

\DeclareMathOperator{\tr}{tr}



\pagestyle{fancy}

\begin{document}

\bibliographystyle{plain}

\else

\fi


\section{Wrapping Number Theorem}
We have found a lower bound for $\tr(K)$ using the winding number, but we still would like to find a stronger bound using the wrapping number. One reason for this is that if the winding number $n = 0$, then Theorem \ref{thm_bound_with_winding_number} does not give anything nontrivial; however, the wrapping number $m$ is always positive, so any bound using it will be nontrivial. The whole proof of Theorem \ref{thm_bound_with_winding_number} works well with the wrapping number except Lemma \ref{lem_1i_1}. This occurs because the proof of Lemma \ref{lem_1i_1} uses the homology interpretation of winding number, and there is no analogous interpretation of the wrapping number. However, we can prove the following key lemma in place of Lemma \ref{lem_1i_1} and conclude our main theorem.

\subsection{Definitions and Lemmas}

\begin{lem}\label{lem_key_lemma}
Let $l'$ be the result of the isotopy of $l$ as in Lemma \ref{lem_li_2}. Suppose $r$ is a regular level of both $h|_{\partial{V}}$ and $h|_{l'}$, and assume that
$$|h^{-1}(r)\cap l'|=s.$$
Let $Q_1,...,Q_s$ be the components of $h^{-1}(r)\cap V$ which intersect $l'$. Then at least $\frac{s}{2}$ of them have exactly one essential boundary component.
\end{lem}

\begin{figure}[H]
    \centering
    \includegraphics[scale=0.25]{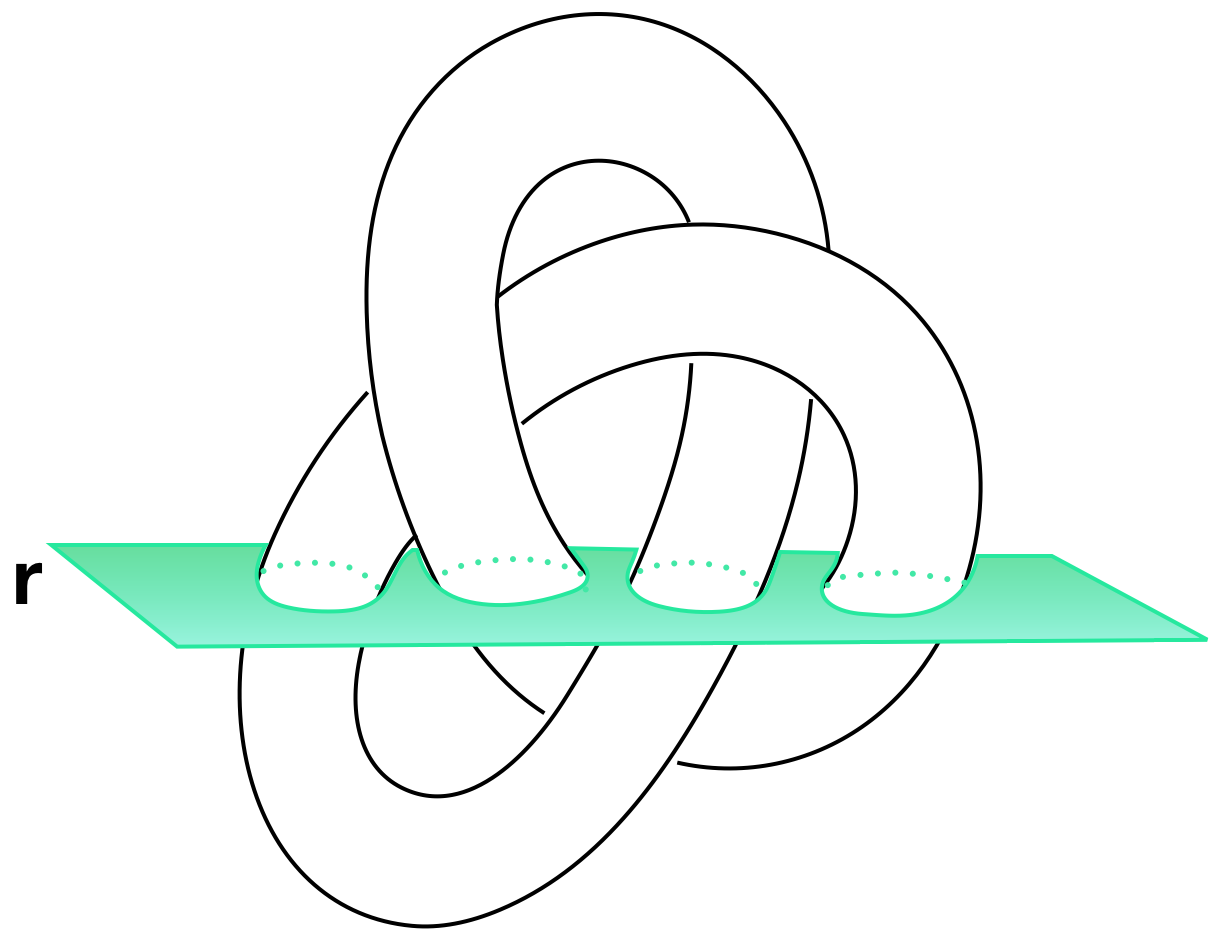}
    \caption{Intersection of a regular level with the solid torus.}
    \label{fig:my_interseciton}
\end{figure}

\begin{thm}\label{thm_trunk_wrapping}
Suppose $K$ is the satellite knot with companion $J$ and the wrapping number is $m$. If $J$ is not the trivial knot, then we have: $\tr(K) > \frac{1}{2}m \cdot \tr(J)$
\end{thm}

Note in this section we will use the lowercase letter $j$ for indices, rather than for the companion knot.

\begin{proof}[Proof of Theorem \ref{thm_trunk_wrapping} by Lemma \ref{lem_key_lemma}.]
Let $k$ be a knot where $\tr(k) = \tr(K)$ and we have the corresponding companion knot and its tubular neighborhood $V$. We construct a loop $l\subset V$ just like in Theorem \ref{thm_bound_with_winding_number}. We can also isotope $l$ into $l'$ as in Lemma \ref{lem_li_2}. Pick a regular level $r$ so that $|h^{-1}(r)\cap l'|= \tr(l')$. We can also look at the components of $h^{-1}(r)\cap V$. Note that each component intersects $l'$ at most once. By Lemma \ref{lem_key_lemma}, among those pieces which do intersect $l'$, more than $\frac{1}{2}$ of them have exactly one essential boundary component. From Definition \ref{winding}, each of these pieces intersects the knot $k$ at least $m$ times (the inessential boundaries can be capped off by disks arbitrarily closed to $\partial{V}$ but $k\subset {\rm int}(V)$ so those piece can be viewed as meridian disks when studying their intersection with $k$). So $\tr(K) > \frac{1}{2}m \tr(L)$. From Theorem \ref{thm_zupan}, we have $\tr(L) \geq \tr(J)$, proving Theorem \ref{thm_trunk_wrapping}.
\end{proof}

In Section 4.1, we will introduce more definitions and lemmas to aid in the proof of Lemma \ref{lem_key_lemma}. This reduces the proof of Lemma \ref{lem_key_lemma} to proving Lemma \ref{lem_key_lemma_comb}, which we do in Section 4.2. We also generalize Theorem \ref{thm_trunk_wrapping} in Section 4.2.

Suppose $V\subset S^3$ is a solid torus so that $h|_{\partial{V}}$ is Morse. Suppose $r\in\real$ is a regular level of $h|_{\partial{V}}$ and
$$h^{-1}(r)\cap V=P_1\cup...\cup P_t.$$

The pieces $P_i$ are contained on the regular level $h^{-1}(r)$ which is diffeomorphic to a $2$-sphere. There are some restrictions on the pieces from topological side. With those restrictions, the question can be solved entirely using combinatorics. 

\begin{lem}\label{lem_requirement_1}
Suppose $P_i$ is as above a component of $h^{-1}(r)\cap V$, then $|P_i \cap l'|= 1$ if and only if $P_i$ has an odd number of essential boundary components.
\end{lem}

\begin{proof}
For each $i$, let $\al_{i,j}$ be essential boundary components and $\be_i$ be the collection of inessential circles. Then 
$$\partial{P}_{i}=\al_{i,1}\cup...\cup\al_{i,s_i}\cup\be_i.$$

We have a boundary map:
$$\partial: H_2(V,\partial{V})\ra H_1(\partial{V})$$
and this map can be described explicitly as follows. $H_1(\partial{V})\cong \intg\oplus \intg$, and the two generators are represented by meridians and longitudes (with respect to some framing). Note also $H_2(V,\partial{V})\cong\intg$ so the map $\partial$ is actually:
$$\partial(1)=(1,0).$$

Then $|P_i\cap l'|=1$ means that $[P_i,\partial{P}_i]=\pm1 H_2(V,\partial{V})$. From the definition of boundary map, we have
\beq
(\pm1,0)&=\partial(\pm1)=\partial([P_i,\partial{P}_i])=[\partial{P}_i]=\sum_{j = 1}^{s_i} [\al_{i, j}]
\eeq

Since $\al_{i,j}\cap \al_{i,j'}=\emptyset$, we know that all $\al_{i,j}$, if given the correct orientation, would represent the same class in $H_2(\partial{V})$. So suppose for all $j$, $[\alpha_{i, j}] = \pm(x, y)$. Since some of the $[\alpha_{i, j}]$ cancel each other out because of opposite orientations, we get 
 $$[\partial P_i] = \sum_{j = 1}^{s_i} [\al_{i, j}] = l \cdot (x, y) = (1, 0)$$
for $|l| \leq s_i$. This implies $l x =\pm 1$ so $l = \pm 1$. Since $l \equiv s_i \pmod 2$, we have that $s_i$ is odd as desired. 

On the other hand, if $P_i$ has an odd number of essential boundary components, then by the above argument it must represent a non-trivial element in $H_2(V,\partial{V})$. Since $l'$ is a generator of $H_1(V)$, we know that $P_i\cap l'\neq \emptyset$. From the construction of $l'$, we know that $|P_i \cap l'|= 1$.
\end{proof}

\begin{lem}\label{lem_essential_boundary_contains_another_piece_inside}
Suppose the solid torus $V$, and hence the companion knot $J$, is knotted. Let $r$ be a regular level and suppose $Q_i$ is a component of $h^{-1}(r)\cap V$ which intersects $l'$. Let $\al$ be an essential boundary component of $Q_i.$ Then $\al$ bounds a disk $D\subset h^{-1}(r)$ whose interior is disjoint from $Q_i.$ Further, there exists another component $Q_j \subset D$ of $h^{-1}(r)\cap V$ which also intersects $l'.$
\end{lem}

\bpf
If $D$ does not contain any piece $Q_j$ in its interior, then 
$$D \cap (S^3 \backslash \mathring{V}) = \partial D = \al.$$
Since $\al$ is essential, this means that the complement of $V$ in $S^3$ is compressible. However, this is absurd since $V$ is knotted. 

As a result, we know that there is a component $P_j$ of $h^{-1}(r)\cap V$ so that $P_j\subset D$. If $P_j$ intersects $l'$ then we are done. If not, and if $P_j$ does not contain any essential boundary components, then we can use the above trick again to cap off all boundary components of $P_j$. If $P_j$ has essential boundary components, then by Lemma \ref{lem_requirement_1}, $P_j$ has at least two essential boundary components. Hence there is an essential boundary component $\be$ of $P_j$ so that $\be$ bounds a disk $D'\subset D$ which is disjoint from the interior of $P_j$. We apply the above argument on $D'$ and thus there is another component $P_j'$ of $h^{-1}(r)\cap V$ which is contained in $D'$. We can keep this argument running until we find a piece $Q_j$ which intersects $l'$.


\epf

Now we introduce the combinatorial setup.

\begin{defn}\label{defn_arrangement}
Suppose $A(s)$ is an embedding of $s$ many compact connected surfaces, or what we called {\it pieces}, into $S^2$:
$$P_1\sqcup P_2 ... \sqcup P_s\hookrightarrow S^2,$$
so that all of the following hold:

(1). For each $i$, let the boundary components of $P_i$ be
$$\partial{P_i}=\al_{i,1}\cup...\cup \al_{i,s_i}.$$
Then $s_i$ is either $1$ or at least $3$.

(2). On the sphere $S^2$, each $\al_{i,j}$ bounds a disk $D_{i,j}\subset S^2$ so that $D_{i,j}\cap {\rm int}(P_i)=\emptyset$. We have the following two requirements on $D_{i,j}$:

(i). For each $D_{i,j}$, there exists a piece $P_x$ so that $P_x\subset{\rm int}(D_{i,j})$.

(ii). If a $D_{i,j}$ does not contain any other disks $D_{i',j'}$, then there exists a piece $P_k$ so that $P_k\subset D_{i,j}$ and $s_k=1$.

We call such $A(s)$ an {\it arrangement} (of the surfaces on a sphere). 
\end{defn}

\begin{defn} For an arrangement $A(s)$, let $\lambda(A(s))$ be the number of pieces which have exactly one boundary component.
\end{defn}

\begin{figure}[H]
   \begin{minipage}{0.48\linewidth}
  \centerline{\includegraphics[width=4.5cm]{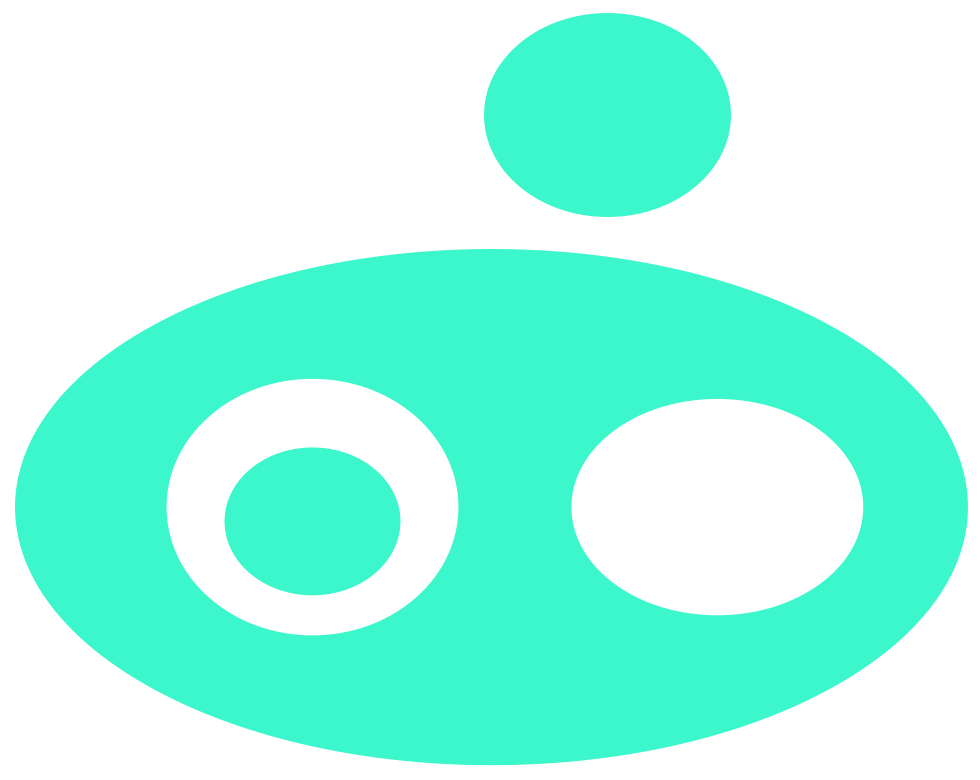}}
  \centerline{Not an arrangement}
\end{minipage}
\hfill
\begin{minipage}{.48\linewidth}
  \centerline{\includegraphics[width=4.5cm]{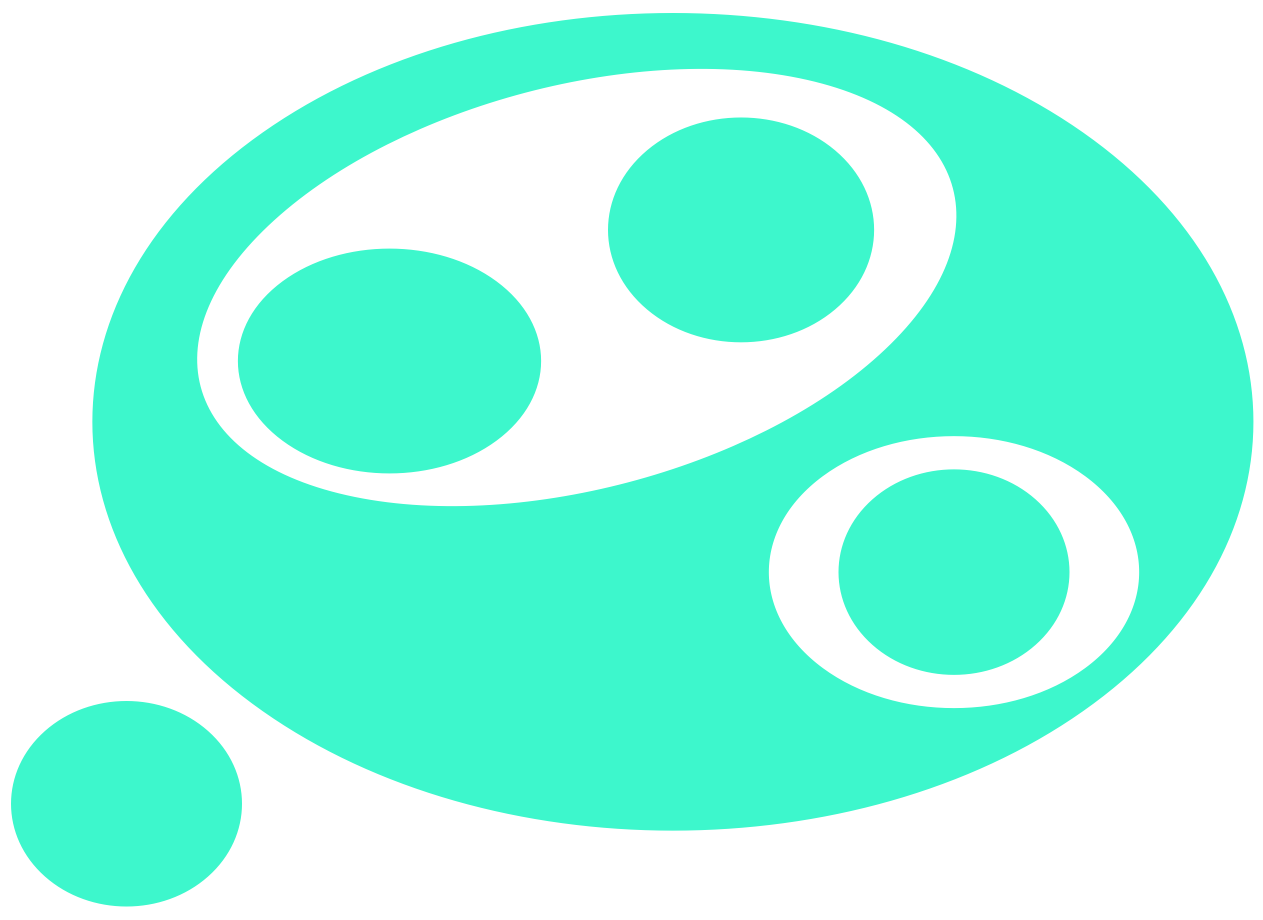}}
  \centerline{Arrangement}
\end{minipage}
   \caption{Examples of piece configurations.}
   \label{fig_arrangement}
\end{figure}

\begin{exmp}
In Figure \ref{fig_arrangement}, we have two pictures. On the left we do not have an arrangement since one innermost piece is missing, violating requirement (2) in Definition \ref{defn_arrangement}. On the right, we have an arrangement where $\lambda(A(5)) = 4$.
\end{exmp}

\begin{defn}\label{defn_interior}
Suppose $P_i$ is a piece and $\be\subset\partial{P}_i$ is a boundary component of $P_i$ which is inessential. Then there is a (unique) disk $D\subset S^2$ so that $D\cap{\rm int}(P_i)=\emptyset$. We call $\be$ {\it pseudo-essential} if $D$ contains another piece $P_j$ in its interior.
\end{defn}

\begin{lem}\label{lem_constructing_an_arrangement}
Under the setup of Lemma \ref{lem_key_lemma}, assume that $t$ is the number of pieces $Q_i$ among $Q_1,...,Q_s$, so that 

(1). $Q_i$ has a unique essential boundary component.

(2). There is a unique pseudo-essential boundary component of $Q_i$ that bounds a disk $D\subset h^{-1}(r)$ whose interior is disjoint from $Q_i.$ Further, there exists a piece $Q_j$ so that $Q_j\subset D$.

Assume without loss of generality that $Q_{s-t+1},...,Q_{s}$ are the $t$ pieces which satisfy conditions (1) and (2) above. Then there is an arrangement 
$$A(s-t)=P_1\cup...\cup P_{s-t}\hookrightarrow h^{-1}(r)\cong S^2,$$
so that for all $i=1,...,s-t$, we have
$$Q_i\subset P_i$$
and $P_i$ has exactly one boundary component only if $Q_i$ has exactly one essential boundary component.
\end{lem}

\begin{proof}
Recall that $Q_1,..,Q_s$ are all components of $h^{-1}(r)\cap V$ which intersect $l'$ and $Q_{s-t+1},...,Q_s$ are the components which satisfy conditions (1) and (2) in the hypothesis of the lemma.  We construct $P_i$ for $i=1,...,s-t$ as follows.

For any $Q_i$, let $\al\subset\partial Q_i$ be an inessential boundary component. Then $\al$ bounds a disk $D\subset h^{-1}(r)\cong S^2$ whose interior is disjoint from $Q_i$. If $D$ does not contain any other pieces $Q_j$ for $i \leq j \leq s- t$, we glue $D$ to $Q_i$. Perform this operation for any possible inessential boundary components of $Q_i$ and let the result be $P_i$. From the construction it is clear that the following claim holds:

{\bf Claim 1}. The number of essential boundary components of $Q_i$ is no larger than the number of boundary components of $P_i$.

As a direct result, if $P_i$ has one boundary component, it follows that $Q_i$ must have only one essential boundary component.

Now we prove that $P_1\cup...\cup P_{s-t} \subset h^{-1}(r) \cong S^2$ is an arrangement as defined in Definition \ref{defn_arrangement}. Note from Lemma \ref{lem_requirement_1} that each $Q_i$ must have an odd number of essential boundaries. By Claim 1, if $Q_i$ has three essential boundary components, then $P_i$ has at least three boundary components. If $Q_i$ has one essential boundary component, then from the hypothesis, since $i\leq s-t$, there are at least two pseudo-essential boundary components $\al'_{i,1}$ and $\al'_{i,2}$ so that:

(i). Each $\al'_{i,j}$ bounds a disk $D'_{i,j}\subset h^{-1}(r)$ whose interior is disjoint from $Q_j$.

(ii). Each $D'_{i,j}$ contains some other pieces $Q_{z_j}$ for $1 \leq z_j \leq s- t$.

Hence, by construction, $\al'_{i,1}$ and $\al'_{i,2}$ survive in $P_i$. So together with the essential boundary component of $Q_i$, $P_i$ has at least three boundary components. Thus condition (1) of Definition \ref{defn_arrangement} is satisfied. 

To prove that condition (2) of Definition \ref{defn_arrangement} is also satisfied, we proceed as follows. For each $i=1,...,s-t$, suppose
$$\partial P_{i}=\al_{i,1}\cup...\al_{i,d}.$$
Each circle $\al_{i,j}$ bounds a disk $D_{i,j}\subset h^{-1}(r)$ whose interior is disjoint from $P_j.$
If $\al_{i,j}$ is an essential boundary component of $Q_i$, then we can apply Lemma \ref{lem_essential_boundary_contains_another_piece_inside} to get that there is a piece $Q_z\subset D_{i,j}$. If $z\leq s-t$, then this means that $P_z\subset D_{i,j}$. If $z > s-t$, then from either the condition (2) in the hypothesis of the current lemma, or from Lemma \ref{lem_essential_boundary_contains_another_piece_inside}, we know that $Q_z$ has a boundary component $\al$ which bounds a disk $D\subset D_{i,j}$ so that there is another piece $Q_{z'}\subset D\subset D_{i,j}$. Hence we can look at $Q_{z'}$ and apply the same innermost argument once more. Since there in total there finitely many $Q_i's$, we will eventually find $Q_{z''}\subset D_{i,j}$ with $z''\leq s-t$. This proves (i) in condition (2) of Definition \ref{defn_arrangement}. The term (ii) of condition (2) actually follows from the same type of innermost argument.
\end{proof}

\begin{lem}\label{lem_key_lemma_comb}
For any arrangement $A(s)$, we have
$$\lambda(A(s))>\frac{s}{2}.$$
\end{lem}

\begin{proof}[Proof of Lemma \ref{lem_key_lemma} by Lemma \ref{lem_key_lemma_comb}]
Under the hypothesis of Lemma \ref{lem_key_lemma}, assume as above that $Q_{s-t+1},...,Q_{s}$ are the pieces satisfying condition (1) and (2) in Lemma \ref{lem_constructing_an_arrangement}. Then by that lemma, we can find an arrangement $A(s-t)$. By Lemma \ref{lem_key_lemma_comb} at least $\frac{1}{2}(s-t)$ pieces of $P_i$ has exactly one boundary component. Then from the conclusion of Lemma \ref{lem_constructing_an_arrangement}, we know that among $Q_1,...,Q_{s-t}$, there are at least $\frac{1}{2}(s-t)$ of them having exactly one essential boundary components. Note all of $Q_{s-t+1},...,Q_{s}$ have exactly one essential boundary component so we are done.
\end{proof}

To finally conclude the proof of Theorem \ref{thm_trunk_wrapping}, we still need to prove Lemma \ref{lem_key_lemma_comb}. In the following subsection, we define a new combinatorial invariant that generalizes the wrapping number. We then prove Theorem \ref{thm_lambda_a}, which implies Lemma \ref{lem_key_lemma_comb} as a corollary.

\subsection{$\lambda(a)$ and $\mu(a)$}

Theorem \ref{thm_trunk_wrapping} establishes a lower bound for the trunk of a satellite knot in terms of the trunk of the companion knot and the wrapping number of the satellite knot. We know from Definition \ref{winding} that the wrapping number is defined to be the minimal geometric intersection number of a meridian disk with the satellite knot. Recall that a piece is a meridian disk if and only if it has exactly $1$ essential boundary. In this section, we extend previous results to obtain a lower bound on the trunk of a satellite knot in terms of the minimal number of geometric intersections of pieces with more than one essential boundary with the satellite knot.

\begin{defn} $S(a)$ is the set of all connected planar surfaces $S \subset V$ such that $\partial S \subset \partial V$, $S$ represents a generator of $H_2(V, \partial V)$ and $S$ has no more than $a$ essential boundaries.

\end{defn}

\begin{defn} Let $\mu(a) = \displaystyle{\min_{S \in S(a)}} |S \cap k|$. 

\end{defn}

By definition, $\mu(a + 1) \leq \mu(a)$. Noting this, we let $$\mu = \lim_{a\to\infty} \mu(a).$$

\begin{defn}\label{lambda} Let $\lambda(a)$ be the largest possible value such that for any satellite knot $k$ with companion $j$ we have $\tr(K) \geq \lambda(a) \cdot \mu(a) \cdot \tr(J)$.

\end{defn}

\begin{rem}\label{winding_relationship}
Note we always have $m=\mu(1) \geq \mu \geq n$, where $n$ is the winding number and $m$ is the wrapping number. There indeed exist cases when $\mu=m$ or $\mu=n$. For instance, when the satellite knot $K$ is the Whitehead double knot, by the work of Li and Guo in \cite{li2018width}, we have that $\mu(a) = 2$ for any $a$. Further, we have that $\tr(K) = 2\tr(J)$, making $\lambda(a) = 1$ in this case. Since we can always form the Whitehead double of any knot, we know that if $\lambda(a) > 1$, then the Whitehead double would be a counterexample. Therefore, we have $\lambda(a) \leq 1$.
\end{rem}

\begin{prop}\label{bounds} $\tr(K) \geq \frac{1}{2}(m + \mu) \cdot \tr(J) \geq \frac{1}{2}(m + n) \cdot \tr(J), \mu \cdot \tr(J)$. 

\end{prop}

Note that Remark \ref{winding_relationship} and Proposition \ref{bounds} imply Theorem \ref{thm_bound_with_winding_number}.

\begin{proof}
Let $k$ be a knot such that $\tr(k) = \tr(K)$ and we have the corresponding companion knot and its tubular neighborhood $V$. We construct a loop $l\subset V$ just like in Theorem \ref{thm_bound_with_winding_number}. We can also isotope $l$ into $l'$ as in Lemma \ref{lem_li_2}. Pick a regular level $r$ for which $|h^{-1}(r)\cap l'|= \tr(l')$. We can also look at the components of $h^{-1}(r)\cap V$. 
The intersection of a regular level $h^{-1}(r)$ with the solid torus is a set of horizontal pieces, each with an odd number of essential boundaries. Let $b_a$ denote the proportion of total pieces with exactly $a$ boundaries. It is obvious that:

$$\sum_{a = 1}^{\infty} b_a = 1.$$


We note that by Lemma \ref{lem_key_lemma}, $b_1 > \frac{1}{2}$. Then we have:
\beq
\tr(K) &\geq \tr(l')\sum_{a = 1}^{\infty} b_a\mu(a)\\
&\geq \tr(l') \left(\frac{1}{2}m + \sum_{a = 3}^{\infty} b_a\mu(a) \right) \geq \tr(l') \left(\frac{1}{2}m + \mu \cdot \sum_{a = 3}^{\infty} b_a \right).  
\eeq


Recalling that $m \geq \mu \geq n$ and $\tr(l') \geq \tr(L) \geq \tr(J)$ we get:

\beq
\tr(K) &\geq \tr(l') \left(\frac{1}{2}m + \mu \cdot \sum_{a = 3}^{\infty} b_a \right) \geq \tr(J) \cdot \frac{\mu + m}{2}\\
&\geq \frac{1}{2}(m + n)\tr(J), \mu \cdot \tr(J).
\eeq

\end{proof}


To bound $\lambda(a)$, we will slightly alter Definition \ref{defn_arrangement} for arrangement: 

\begin{defn}\label{defn_arrangement_changed}
Suppose $a$ is an odd positive integer number. An {\it arrangement} $A(s)$ is an embedding of surfaces as in Definition \ref{defn_arrangement}, but with condition (1) replaced by the following:

(1)'. For each $i$, let the boundary components of $P_i$ be
$$\partial{P_i}=\al_{i,1}\cup...\cup \al_{i,s_i}.$$
Then either $s_i\leq a$ or $s_i\geq a+2$.
\end{defn}



\begin{lem}\label{lem_two_moves}
Any arrangement as defined in Definition \ref{defn_arrangement_changed} can be constructed from two disks using a sequence of the following two types of ``moves". 

Move 1: adding a new disk to the arrangement. 

Move 2: take $s\in\intg$ so that $0\leq s\neq a+1$. Then replacing a disk with a piece with $s$ boundary components and also adding one disk inside each of the newly created boundary components. 
\end{lem}

Note that each step is reversible by definition, and both performing the move and reversing the move in an arrangement results in another arrangement. Also, the initial state consisting of two disks is an arrangement.

\begin{proof}[Proof of Lemma \ref{lem_two_moves}]

Given any arrangement $A(s)$ (pieces on a $2$-sphere), we repeat two procedures until it is not possible to continue: 

Procedure 1: perform the reverse of Move 1 until each boundary component contains at most one disk. 

Procedure 2: perform the reverse of Move 2 on pieces with more than 1 boundary component that do not contain another piece with more than 1 boundary. 

If there exists a piece with more than 1 boundary component, then by the innermost argument, there exists a piece with more than 1 boundary component such that each boundary component contains only disks within its interior (see Definition \ref{defn_interior} for interior). After performing the first procedure to the fullest extent, the second procedure can be performed to remove the ``innermost" piece with more than one boundary component. Thus, the only possible arrangement where these procedures cannot be performed does not have any piece with more than one boundary component, meaning it only has disks. After reaching such an arrangement, perform procedure 1 to get two disks. 

To build the original arrangement from these two remaining disks, reverse each move made in the deconstruction process.
\end{proof}

\begin{thm}\label{thm_lambda_a} Suppose $a$ is odd. Then we have $\lambda(a) > \frac{a}{a+1}$.

\end{thm}

\begin{proof}
We use an inductive argument. Let $A(s_t)$ be the arrangement after $t$ moves have been performed on an arrangement $A(s)$. Let $x_t$ be the number of pieces with at most $a$ boundaries, and let $y_t$ be the total number of pieces in $A(s_t)$. From Definition \ref{lambda} and the same argument as in the proof of Theorem \ref{thm_trunk_wrapping}, $\lambda(a)$ is greater than or equal to the minimum possible value of $\frac{x_t}{y_t}$. Note that we start with $A(s_0)=A(2)$, which is the arrangement with two disks, and in this case, $x_0=2$ and $y_0=2$.

Performing Move 1 results in $x_{t+1} = x_t + 1$ and $y_{t+1} = y_t + 1$ so $\frac{x_{t+1}}{y_{t+1}} \geq \frac{x_t}{y_t}$. 

Performing Move 2 with a piece of $c$ boundaries such that $c \leq a$ results in $x_{t+1} = x_t + (c-1)$ and $y_{t+1} = y_t + (c-1)$ so $\frac{x_{t+1}}{y_{t+1}} \geq \frac{x_t}{y_t}$.

Performing Move 2 with a piece of $c$ boundaries such that $c > a$ results in $x_{t+1} = x_t + (c-2)$ and $y_{t+1} = y_t + (c-1)$. Note that by condition (1)' in definition \ref{defn_arrangement_changed}, $c\geq a+2$. A smaller number of essential boundaries in the piece added always causes the biggest decrease in $\frac{x_t}{y_t}$, as $\frac{x_t+d-2}{y_t+c-1} > \frac{x_t+d+2-2}{y_t+d+2-1}$ since $x_t \leq y_t$ by definition. Thus, the arrangement with minimal $\frac{x_t}{y_t}$ is constructed by repeatedly performing Move 2, converting disks to pieces with $a+2$ boundaries: 

$$\lim_{t\to\infty} \frac{x_t}{y_t} = \lim_{t\to\infty} \frac{2+ t \cdot a}{2+t \cdot (a+1)} = \frac{a}{a+1}.$$

Thus, $\lambda(a)>\frac{a}{a+1}$. Note that when $a = 1,$ we get $\lambda(a) > \frac{1}{2},$ proving Lemma \ref{lem_key_lemma_comb}.
\end{proof}


\ifx\allfiles\undefined

\newpage

\bibliography{Index}

\end{document}

\fi


\bibliography{Index}

\begin{thebibliography}{10}

\bibitem{blair2013width}
Ryan {Blair} and Maggy {Tomova}.
\newblock Width is not additive.
\newblock {\em Geometry and Topology}, 17(1):93--156, 2013.

\bibitem{davies2017natural}
Derek {Davies} and Alexander {Zupan}.
\newblock Natural properties of the trunk of a knot.
\newblock {\em Journal of Knot Theory and Its Ramifications}, 26(12):1750080,
  2017.

\bibitem{gabai1987foliations}
David Gabai et~al.
\newblock Foliations and the topology of 3-manifolds. iii.
\newblock {\em Journal of Differential Geometry}, 26(3):479--536, 1987.

\bibitem{gordon1989knots}
C.~McA. {Gordon} and J.~{Luecke}.
\newblock Knots are determined by their complements.
\newblock {\em Journal of the American Mathematical Society}, 2(2):371--415,
  1989.

\bibitem{guo2018width}
Qilong Guo and Zhenkun Li.
\newblock Width of a satellite knot and its companion.
\newblock {\em Algebraic and Geometric Topology}, 18(1):1--13, 2018.

\bibitem{li2018width}
Zhenkun {Li} and Qilong {Guo}.
\newblock Width of the whitehead double of a nontrivial knot.
\newblock {\em arXiv preprint arXiv:1803.02480}, 2018.

\bibitem{ozawa2010waist}
Makoto {Ozawa}.
\newblock Waist and trunk of knots.
\newblock {\em Geometriae Dedicata}, 149(1):85--94, 2010.

\bibitem{rieck2002thin}
Yo'av {Rieck} and Eric {Sedgwick}.
\newblock Thin position for a connected sum of small knots.
\newblock {\em Algebraic and Geometric Topology}, 2(1):297--309, 2002.

\bibitem{rolfsen1976knots}
Dale Rolfsen.
\newblock {\em Knots and links}, volume 346.
\newblock American Mathematical Soc., 1976.

\bibitem{scharlemann20063}
Martin {Scharlemann} and Jennifer {Schultens}.
\newblock 3-manifolds with planar presentations and the width of satellite
  knots.
\newblock {\em Transactions of the American Mathematical Society},
  358(9):3781--3805, 2006.

\bibitem{schubert1954uber}
Horst {Schubert}.
\newblock \"uber eine numerische knoteninvariante.
\newblock {\em Mathematische Zeitschrift}, 61(1):245--288, 1954.

\bibitem{schultens2003additivity}
Jennifer {Schultens}.
\newblock Additivity of bridge numbers of knots.
\newblock {\em Mathematical Proceedings of the Cambridge Philosophical
  Society}, 135(3):539--544, 2003.

\bibitem{thompson1994position}
Abigail Thompson.
\newblock Thin position and the recognition problem for $s^3$.
\newblock {\em Mathematical Research Letters}, 1(5):613--630, 1994.

\bibitem{wu2006position}
Ying-Qing Wu.
\newblock Thin position and essential planar surfaces.
\newblock {\em Proc. Amer. Math. Soc. 132}, 132:3417--3421, 2004.

\bibitem{alex2010lower}
Alexander Zupan.
\newblock A lower bound on the width of satellite knots.
\newblock {\em Topology Proceedings}, 40:179--188, 2012.

\end{thebibliography}

\end{document}